\documentclass{amsart}

\usepackage{amsmath,amssymb,amsthm,a4wide}
\usepackage{pgfplots}
\usepackage{graphicx,tikz}
\newtheorem{theorem}{Theorem}

\newtheorem*{proposition}{Proposition}
\newtheorem*{corollary}{Corollary}

\newtheorem*{lemma}{Lemma}

\theoremstyle{definition}

\theoremstyle{remark}

\begin{document}

\title[]{ Poissonian Pair Correlation and Discrepancy}

\author[]{Stefan Steinerberger}
\address{Department of Mathematics, Yale University, New Haven, CT 06511, USA}
\email{stefan.steinerberger@yale.edu}

\begin{abstract} A sequence $(x_n)_{n=1}^{\infty}$ on the torus $\mathbb{T} \cong [0,1]$ is said to exhibit Poissonian pair correlation if the local gaps behave like the gaps of a Poisson random variable, i.e.
 $$ \lim_{N \rightarrow \infty}{ \frac{1}{N} \# \left\{ 1 \leq m \neq n \leq N: |x_m - x_n| \leq \frac{s}{N} \right\}} = 2s \qquad \mbox{almost surely.}$$ 
We show that being close to Poissonian pair correlation for few values of $s$ is enough to deduce global regularity statements: if, for some~$0 < \delta < 1/2$, a set of points $\left\{x_1, \dots, x_N \right\}$ satisfies
$$ \frac{1}{N}\# \left\{ 1 \leq m \neq n \leq N: |x_m - x_n| \leq \frac{s}{N} \right\} \leq (1+\delta)2s \qquad \mbox{for all} \hspace{6pt} 1 \leq s \leq (8/\delta)\sqrt{\log{N}},$$
then the discrepancy $D_N$ of the set satisfies $D_N \lesssim \delta^{1/3} + N^{-1/3}\delta^{-1/2}$. 
We also show that distribution properties are reflected in the global deviation from the Poissonian pair correlation
$$ N^2 D_N^5 \lesssim \frac{2}{N}\int_{0}^{N/2} \left|  \frac{1}{N}\# \left\{ 1 \leq m \neq n \leq N: |x_m - x_n| \leq \frac{s}{N} \right\} - 2s \right|^2 ds \lesssim N^2 D_N^2,$$
where the lower is bound is conditioned on $D_N \gtrsim N^{-1/3}$. The proofs use a connection between exponential sums, the
heat kernel on $\mathbb{T}$ and spatial localization. Exponential sum estimates are obtained as a byproduct.
We also describe a connection to diaphony and several open problems.

\end{abstract}

\maketitle

\vspace{10pt}

\section{Introduction and main results}

\subsection{Pair correlation.}
 Let $(x_n)_{n=1}^{\infty}$ be a sequence on the one-dimensional Torus $\mathbb{T} \cong [0,1]$. A natural object of interest is the behavior of gaps between the first $N$ elements on a local scale. If the sequence is comprised of i.i.d. uniformly distributed random variables, then, for all $s>0$, 
 $$ \lim_{N \rightarrow \infty}{ \frac{1}{N} \# \left\{ 1 \leq m \neq n \leq N: |x_m - x_n| \leq \frac{s}{N} \right\}} = 2s \qquad \mbox{almost surely.}$$ 
 Whenever a deterministic sequence $(x_n)_{n=1}^{\infty}$ has the same property, we say it has Poissonian pair correlation. This notion has been intensively investigated, see e.g. \cite{aist1, aist, heath, rud1, rud2, rud3, walker}. While it is well known that there are many such sequences in a metric sense, there is currently no explicit example known (see Pirsic \& Stockinger \cite{pir}). The existing literature has mostly been concerned with whether a deterministic object exhibits Poissonian pair correlation and there are relatively few results about the notion itself. Only recently  Aistleitner, Lachmann \& Pausinger \cite{aist} and Grepstad \& Larcher \cite{grep} independently established that sequences with Poissonian pair correlation are uniformly distributed on $[0,1]$ (see \cite{stein} for another proof).  We believe that there are many interesting open problems regarding this notion and discuss some of them in this paper.

\subsection{A local result.} The first result shows that being close to Poissonian pair correlation for a small range of values of $s$ can be enough to conclude
global regularity results. We quantify regularity using the discrepancy $D_N$ of a finite point set, defined in the usual manner as the maximal deviation between empirical and uniform distribution
$$ D_N = \sup_{I \subset \mathbb{T}} \left| \frac{\#\left\{x_1, x_2, \dots, x_N\right\} \cap I}{N} - |I|\right|,$$
where the supremum ranges over all intervals $I \subset \mathbb{T}$.

\begin{theorem}
Let $\left\{x_1, \dots, x_N \right\} \subset [0,1]$ and $0 < \delta \ll 1$ such that
$$ \frac{1}{N} \# \left\{ 1 \leq m \neq n \leq N: |x_m - x_n| \leq \frac{s}{N} \right\} \leq (1+\delta)2s \qquad \mbox{for all} \hspace{6pt} 1 \leq s \leq (8/\delta)\sqrt{\log{N}},$$
then the discrepancy of the set satisfies $D_N \lesssim  \delta^{1/3} + \delta^{-1/2}N^{-1/3}$.
\end{theorem}

This should be compared to a result of Grepstad \& Larcher \cite{grep} that being $\delta-$close to Poissonian pair correlation for $s \in \left\{1, \dots, \delta^{-5} \right\}$ implies $D_N \lesssim \delta$. This result and Theorem 1 are clearly of the same flavor but cover somewhat different scaling regimes -- we have no reason to assume that these results are optimal.
There should be many other interesting results along these lines.

\begin{quote} \textbf{Open Problem} (Global regularity via local pair correlation statistics)\textbf{.} What is the smallest range of values of $s$ for which one needs to require approximate Poissonian pair correlation statistics to ensure some regularity of the distribution?
\end{quote}

The proof of Theorem 1 is based on the use of Fourier analysis to obtain an exponential sum estimate: we show that for $0 < \delta \ll 1$ the assumptions of Theorem 1 imply
$$\sum_{k \neq 0 \atop |k| \leq \delta^{3/2} N  }{ \left| \sum_{n=1}^{N}{ e^{2 \pi i  kx_n }}\right|^2}  \lesssim  \delta N^2.$$
It is instructive to study the case of randomly chosen points: then each of the $\sim \delta^{3/2}N$ squared
exponential sums are of expected size $\sim N$ and the expression would be $\sim \delta^{3/2} N^2$. It is not clear to us whether such a bound holds or whether the assumptions in Theorem 1 allow for point sets that are substantually different from randomly chosen points.

\subsection{A global result.}
We show that well-distributed sequences have pair correlation globally close to Poissonian.
For somewhat irregular sets, $D_N \gtrsim N^{-1/3}$, the converse direction also holds.

\begin{theorem} Let $\left\{x_1, \dots, x_N \right\} \subset [0,1]$. Then
$$\frac{2}{N}\int_{0}^{N/2} \left| \frac{1}{N}  \# \left\{ 1 \leq m \neq n \leq N: |x_m - x_n| \leq \frac{s}{N} \right\} - 2s \right|^2 ds \lesssim N^2 D_N^2.$$
Moreover, if $D_N \gtrsim N^{-1/3}$, then
$$\frac{2}{N}\int_{0}^{N/2} \left| \frac{1}{N}  \# \left\{ 1 \leq m \neq n \leq N: |x_m - x_n| \leq \frac{s}{N} \right\} - 2s \right|^2 ds \gtrsim N^2 D_N^5.$$
\end{theorem}

The statement is sharp in the regime $D_N \sim 1$ since upper and lower bound match but it is very clearly not sharp anywhere else. 
In particular, it would be quite nice to see whether one could possibly obtain results of such a flavor for a more restricted range of values of $s$.

\begin{quote}
\textbf{Open problem.} Can Theorem 2 be improved/sharpened/localized?
\end{quote}

\subsection{Concluding Remarks.} 
The proof of Theorem 2 makes use of LeVeque's upper bound \cite{leveque} derived from the Erd\H{o}s-Turan
inequality; the arising exponential sum is sometimes called the diaphony \cite{zinterhof} 
$$ F_N := \left( 2\sum_{k=1}^{\infty}\frac{1}{k^2} \left| \sum_{n=1}^{N}{ e^{2\pi i k x_n}} \right|^2 \right)^{1/2} = \left( \frac{\pi^2 }{2N^2} \sum_{m,n=1}^{N}{\left( (1-2\left\{x_m - x_n\right\})^2 - \frac{1}{3} \right)}\right)^{1/2}.$$
A byproduct of our proof of Theorem 2 is the following Corollary.

\begin{corollary} Let $\left\{x_1, \dots, x_N \right\} \subset [0,1]$. Then
$$A = \frac{2}{N}\int_{0}^{N/2} \left| \frac{1}{N}  \# \left\{ 1 \leq m \neq n \leq N: |x_m - x_n| \leq \frac{s}{N} \right\} - 2s \right|^2 ds$$
is bounded from above by
$$A \leq \pi^{-2} F_{N^2}\left(\left\{x_n - x_m: 1 \leq n,m \leq N\right\}\right)^2+ 1.$$
\end{corollary}
Since there many results dealing with diaphony of deterministic sequence (see for example \cite{dia3, dia2, dia1}), this Corollary could suggest that there might be some hope of getting refined results for the pair correlation of deterministic sequences.
There is another series of results that seem connected in spirit: given a set $\left\{x_1, \dots, x_N\right\} \subset \mathbb{T}$, we may consider the difference set $\left\{x_i - x_j: 1 \leq i,j \leq N\right\}$
and ask how the discrepancy $D_N$ of the set relates to the discrepancy of the difference set $D_{N^2}$. Improving earlier results by Vinogradov \cite{vin} and van der Corput \& Pisot \cite{vdc}, Cassels \cite{cas} showed
$$D_N \lesssim \sqrt{D_{N^2}} ( 1+ |\log{D_{N^2}}|).$$
Motivated by this result, we quicky note another approach that follows rather quickly from the Erd\H{o}s-Turan inequality but may prove useful for such problems or even be of independent interest.

\begin{proposition} Let $\left\{x_1, \dots, x_N\right\} \subset \mathbb{T}$. There is a discrepancy bound
$$D_N  \lesssim \frac{\sqrt{\log{N}}}{N} \sum_{m,n = 1}^{N}{  \min\left\{ \log{N}, \log{\left(\frac{1}{4 \sin{\left(\pi (x_m - x_n)\right)}^2}\right)}\right\}}.$$
\end{proposition}
We observe that
$$ \int_{0}^{1/2}{  \log{\left(  \frac{1}{4 \sin^2{\pi x}}  \right)} dx} = 0,$$
which indicates that there should be cancellation in the sum if the set of points has a pair correlation close to Poissonian. While this approach might not yield
localized estimates it could conceivably lead to results along the lines of Theorem 2. 
Finally, unconnected to these other results, we note two other curious by-product of the proof of Theorem 2. The first is another proof that $\zeta(2) = \pi^2/6$ (given after the proof). The second implication is an exponential sum estimate.
\begin{corollary} For all $\left\{x_1, \dots, x_N\right\} \subset \mathbb{T}$
$$ N\sum_{ k = 0 \atop k~{\tiny \mbox{odd}}}^{\infty}   \frac{8}{k^2 }\left| \sum_{n=1}^{N}{ e^{2 \pi i  k x_n }}\right|^2 \leq  \sum_{k =1}^{\infty}  \frac{2}{k^2} \left| \sum_{n=1}^{N}{ e^{2 \pi i  k x_n }}\right|^4  + \pi^2 N^2.$$
\end{corollary}
One interpretation of that inequality is that any finite set of points
cannot only be irregular with respect to odd frequencies. It could be interesting to see whether this inequality is part of a larger family of inequalities, at least visually it seems to have a certain interpolatory flavor.
We conclude by remarking that a weaker notion was already introduced in \cite{stein}, where it was shown that if $(x_n)_{n=1}^{\infty}$ is a sequence on $\mathbb{T}$, $0 < \alpha < 1$ and for all $s > 0$
  $$ \lim_{N \rightarrow \infty}{ \frac{1}{N^{2-\alpha}} \# \left\{ 1 \leq m \neq n \leq N: |x_m - x_n| \leq \frac{s}{N^{\alpha}} \right\}} = 2s \qquad \mbox{a.s.},$$ 
 then the sequence $(x_n)$ is uniformly distributed. We note that this interpolates between Poissonian pair correlation ($\alpha = 1$) and a classical notion of uniform distribution ($\alpha = 0$).
\begin{quote} \textbf{Open problem.} Do 'most' sequences satisfy this property for some $0 <  \alpha < 1$?
\end{quote}
It seems conceivable that $D_N \lesssim N^{-\beta}$ would imply the property for all $\alpha < \beta$. However, there are other natural notions that could be of interest, we specifically
mention conditions like
$$\int_{s-\frac{1}{2}}^{s+\frac{1}{2}}  \frac{1}{N}  \# \left\{ 1 \leq m \neq n \leq N: |x_m - x_n| \leq \frac{t}{N} \right\} dt \rightarrow_{N \rightarrow \infty} 2s \quad \mbox{for}~s \geq 1/2$$
or, for $s \gg 1$ and $u = o(s)$,
$$ \frac{1}{2u}\int_{s-u}^{s+u}  \frac{1}{N}  \# \left\{ 1 \leq m \neq n \leq N: |x_m - x_n| \leq \frac{t}{N} \right\} dt \rightarrow_{N \rightarrow \infty} 2s + o(u).$$

\section{Proof of Theorem 1}

\subsection{Preliminaries.} We will use the Jacobi $\theta-$function given by
$$ \theta_t(x) = \sum_{k \in \mathbb{Z}}{ e^{-4 \pi^2 k^2 t}e^{2 \pi i k x}} = 1 + 2\sum_{k=1}^{\infty}{e^{-4 \pi^2 k^2 t} \cos{2 \pi k t}}.$$
Basic properties are $\theta_t(x) \geq 0$ and
$$ \int_{\mathbb{T}}{\theta_t(x) dx} = 1.$$
We will use it as a tool that allows us to localize functions: convolution with $\theta_t$ is easy to compute since its Fourier series is explicit. Simultaneously, convolution has little effect on the function since $\theta_t(x)$ is
highly localized: $\theta_t(x)$ acts as the heat kernel on $\mathbb{T}$ and thus, for $t$ small, is well-approximated by the Euclidean heat kernel
$$ \theta_t(x) \sim \frac{1}{\sqrt{4 \pi t}} e^{-\frac{|x|^2}{4t}}.$$
There are various ways of making this notion precise, one of them being that the heat kernel $k_t$ on $\mathbb{R}$ satisfies
$$ k_t(x) =  \frac{1}{\sqrt{4 \pi t}} e^{-\frac{|x|^2}{4t}} \qquad \mbox{and} \qquad \theta_t(x) = \sum_{k \in \mathbb{Z}}{k_t(x+k)}.$$

The second ingredient that we need is a fairly basic rearrangement statement: its underlying idea is far from novel but this particular case
may not have been stated before (though it can be proved in the usual completely standard manner).
\begin{lemma} Let $f:[0, \infty] \rightarrow \mathbb{R}_{\geq 0}$ be a strictly monotonically decreasing function and suppose that the finite measure $\mu$ on $[0, \infty]$ satisfies $\mu\left(\left[0, x\right]\right) \leq \phi(x)$ for all $\alpha < x < \beta$ for some $\phi \in C^1$. Then
$$ \int_{0}^{\infty}{f d \mu} \leq f(0) \phi(\alpha) + \int_{\alpha}^{\beta}{f(x) \phi'(x) dx} + f(\beta)\mu\left( \mathbb{R}_{\geq 0} \setminus [0, \beta]\right).$$
\end{lemma}
The proof is an elementary rearrangement argument, see e.g. Lieb \& Loss \cite{lieb}, and is left to the reader. Indeed, it is not difficult to see that
the right-hand side is sharp and the extremal measure $\mu$ can be characterized: it has point mass $\phi(\alpha)$ in 0, the absolutely continuous density $\phi'(x) dx$ 
on $[\alpha, \beta]$ and another point mass at $\beta$. We will not use the characterization and, when applying the result, replace the last term by the larger quantity $ f(\beta)\mu\left( \mathbb{R}_{\geq 0}\right)$.

\subsection{Proof of the Theorem}
\begin{proof} 
For $t$ small, one summand dominates the remaining expression. We start the argument by using an idea from \cite{stein, stein1}: for arbitrary $X> 0$
\begin{align*}
\sum_{|k| \leq \delta^{3/2} N}{ \left| \sum_{n=1}^{N}{ e^{2 \pi i  k x_n }}\right|^2}  &\leq 
e^{4\pi^2 \delta} \sum_{k \in \mathbb{Z}}{e^{-4 \pi^2 k^2/ (\delta N)^2}  \left| \sum_{n=1}^{N}{ e^{2 \pi i  k x_n }}\right|^2} \\
&= e^{4\pi^2 \delta}\sum_{k \in \mathbb{Z}}{e^{-4 \pi^2 k^2/ (\delta N)^2}   \sum_{m, n=1}^{N}{ e^{2 \pi i  k(x_m-x_n) }}} \\ 
&=  e^{4\pi^2 \delta}  \sum_{m, n=1}^{N} \sum_{k \in \mathbb{Z}}{e^{-4 \pi^2 k^2/ (\delta N)^2} e^{2 \pi i  k (x_m-x_n)}}\\ 
&= e^{4\pi^2 \delta}\sum_{m,n=1}^{N}{\theta_{(\delta N)^{-2}}(x_m - x_n)}.
 \end{align*}
We introduce the measure (given as the finite sum of Dirac measures)
$$ \mu = \sum_{m,n=1 \atop m \neq n}^{N}{ \delta_{x_m - x_n}}$$
and use it to write
\begin{align*}
\sum_{m,n=1}^{N}{\theta_{(\delta N)^{-2}}(x_m - x_n)} &= N\theta_{(\delta N)^{-2}}(0) +  \int_{\mathbb{T}}^{}{\theta_{(\delta N)^{-2}}(x) d\mu }.
\end{align*}
The function $\theta_{(\delta N)^{-2}}$ is monotonically decaying away from the origin and the measure $\mu$ satisfies 
$$ \mu\left(\left[ -s, s \right]\right) \leq \left(1 + \delta\right) 2 s N^2 \qquad \mbox{for all}~\quad \frac{1}{N} \leq s \leq \frac{8}{\delta}\frac{\sqrt{\log{N}}}{N}.$$
This implies, using the symmetry of $\mu$ and the previous Lemma,
\begin{align*}
  \int_{\mathbb{T}}^{}{\theta_{(\delta N)^{-2}}(x) d\mu } &\leq   \theta_{(\delta N)^{-2}}(0) (1+\delta) 2 N + (1+\delta)2 N^2  \int_{1/N}^{(8/\delta) \sqrt{\log{N}}/N}{ \theta_{(\delta N)^{-2}}(x)  dx}\\
&+ \theta_{(\delta N)^{-2}}\left(\frac{8 \sqrt{\log{N}}}{ \delta N}\right) N^2.
\end{align*}
We observe that 
$$  \int_{1/N}^{(8/\delta) \sqrt{\log{N}}/N}{ \theta_{(\delta N)^{-2}}(x) dx} \leq  \frac{1}{2}\int_{\mathbb{T}}^{}{ \theta_{(\delta N)^{-2}}(x)  dx} = \frac{1}{2}.$$
We also observe that
$$ \theta_{(\delta N)^{-2}}\left(\frac{8 \sqrt{\log{N}}}{\delta N}\right) N^2  \leq \left(1+o(1)\right) \frac{\delta N}{\sqrt{4\pi}} \exp\left( -\frac{\delta^2 N^2}{4} \frac{64}{N^2}\frac{ \log{N}}{\delta^2}\right) N^2 \ll 1.$$
Furthermore
$$ \theta_{(\delta N)^{-2}}(0) \sim (1+o(1)) \frac{\delta N}{\sqrt{4\pi}}$$
where $o(1) \rightarrow 0$ as $N \rightarrow \infty$. Summing up, we obtain
$$\sum_{|k| \leq \delta^{3/2} N  }{ \left| \sum_{n=1}^{N}{ e^{2 \pi i  k x_n }}\right|^2}  \lesssim e^{4\pi^2 \delta} \left(  \delta N^2 +  (1+\delta) N^2 \right).$$
and thus, subtracting the value $N^2$ coming from $k=0$,
$$\sum_{k \neq 0 \atop |k| \leq \delta^{3/2} N  }{ \left| \sum_{n=1}^{N}{ e^{2 \pi i  kx_n }}\right|^2}  \lesssim \delta N^2.$$
We can now employ LeVeque's inequality \cite{leveque} to conclude that
$$
 N D(N)\leq \left( N \sum_{k=1}^{\infty}{ \frac{1}{k^2} \left| \sum_{n=1}^{N}{ e^{2 \pi i k x_n }}\right|^2 } \right)^{1/3} \lesssim  N^{1/3} \left(  \delta N^2 + \sum_{k = \delta^{3/2} N}^{\infty}\frac{1}{k^2} \left| \sum_{n=1}^{N}{ e^{2 \pi i  k x_n }}\right|^2 \right)^{1/3}.$$
We use the trivial estimate $\leq N^2$ on the remaining exponential sum 
$$N^{1/3} \left( \delta N^2 + \sum_{k = \delta^{3/2} N}^{\infty}\frac{1}{k^2} \left| \sum_{n=1}^{N}{ e^{2 \pi i  k x_n }}\right|^2 \right)^{1/3} \leq N^{} \left(\delta +   \sum_{ k = \delta^{3/2} N}^{\infty}\frac{1}{k^2}  \right)^{1/3} \lesssim
 N^{} \left(\delta +   \frac{1}{\delta^{3/2} N}  \right)^{1/3}.$$
It is easily seen that 
$$  \left(\delta +   \frac{1}{\delta^{3/2} N}  \right)^{1/3} \lesssim \begin{cases} \delta^{1/3} \qquad &\mbox{if}~\delta \gtrsim N^{-2/5} \\
N^{-1/3} \delta^{-1/2} \qquad &\mbox{if}~\delta \lesssim N^{-2/5}.
\end{cases}
$$

\end{proof}

It is easy to pinpoint where the argument is lossy: in the absence of more information, we assume that the measure $\mu$ is clustered immediately outside of $s= (8/\delta)\sqrt{\log{N}}$. If this
could be excluded, then further improvements could be obtained from the same argument.

\section{Proof of Theorem 2.}

\begin{proof} We use $\chi$ to denote, as usual, the characteristic function of a set and start by rewriting the problem as (note the change of scale $s/N \rightarrow s$)
\begin{align*}
  \# \left\{ 1 \leq m \neq n \leq N: |x_m - x_n| \leq s \right\} = \left\langle \left(\sum_{i=1}^{N}{\delta_{x_i}}\right) * \chi_{\left[-s,s\right]}, \sum_{i=1}^{N}{\delta_{x_i}}\right\rangle - N.
\end{align*}
Plancherel's theorem implies
$$  \left\langle \left(\sum_{i=1}^{N}{\delta_{x_i}}\right) * \chi_{\left[-s,s\right]}, \sum_{i=1}^{N}{\delta_{x_i}}\right\rangle - N = \sum_{k \in \mathbb{Z}} \widehat \chi_{\left[-s,s\right]}(k)\left| \sum_{n=1}^{N}{ e^{2 \pi i  k x_n }}\right|^2 - N$$
and removing the frequency $k=0$ allows us to rewrite the expression as
$$ \sum_{k \in \mathbb{Z}} \widehat \chi_{\left[-s,s\right]}(k)\left| \sum_{n=1}^{N}{ e^{2 \pi i  k x_n }}\right|^2 - N=  \sum_{k \in \mathbb{Z} \atop k \neq 0} \widehat \chi_{\left[-s,s\right]}(k)\left| \sum_{n=1}^{N}{ e^{2 \pi i  k x_n }}\right|^2 + 2sN^2 - N$$
This implies that the quantity
$$ \frac{A}{2} = \int_{0}^{1/2} \left( \frac{1}{N} \# \left\{ 1 \leq m \neq n \leq N: |x_m - x_n| \leq s \right\} - 2Ns \right)^2 ds$$
can be written as
$$ \frac{B}{2}=  \int_{0}^{1/2}{ \left( \frac{1}{N} \sum_{k \in \mathbb{Z} \atop k \neq 0} \widehat \chi_{\left[-s, s\right]}(k)\left| \sum_{n=1}^{N}{ e^{2 \pi i  k x_n }}\right|^2 - 1\right)^2 ds}.$$
We square the expression and deal with the three terms separately. The first term is
\begin{align*}
 \frac{1}{N^2} \int_{0}^{1/2}   \sum_{k,m \in \mathbb{Z} \atop k \neq 0 \neq m} \widehat \chi_{\left[-s,s\right]}(k) \widehat \chi_{\left[-s,s\right]}(m) \left| \sum_{n=1}^{N}{ e^{2 \pi i  k x_n }}\right|^2 \left| \sum_{n=1}^{N}{ e^{2 \pi i  m x_n }}\right|^2   ds,
\end{align*}
which can be rearranged as
$$ \frac{1}{N^2}  \sum_{k,m \in \mathbb{Z} \atop k \neq 0 \neq m} \left| \sum_{n=1}^{N}{ e^{2 \pi i  k x_n }}\right|^2 \left| \sum_{n=1}^{N}{ e^{2 \pi i  m x_n }}\right|^2  \int_{0}^{1/2} \widehat \chi_{\left[-s,s\right]}(k) \widehat \chi_{\left[-s,s\right]}(m) ds.$$
We quickly compute all arising integrals: for $k \in \mathbb{Z}$ and $k \neq 0$,
$$  \widehat \chi_{\left[-s, s\right]}(k) =    \int_{-1/2}^{1/2}{\chi_{|y| \leq s} e^{-2 \pi i k y} dy} ds = \frac{\sin{(2 k \pi s)}}{k \pi}.$$
It is easy to see that the expression vanishes when integrated over $[0,1/2]$ if $k$ is even (and $k \neq 0$). If $k$ is odd, then
$$ \int_{0}^{1/2}{ \frac{\sin{(2 k \pi s)}}{k \pi} ds} =  \int_{\frac{k-1}{2k}}^{1/2}{ \frac{\sin{(2 k \pi s)}}{k \pi} ds} = \frac{1}{k^2 \pi^2}.$$
Moreover, for $k,m \in \mathbb{Z}$, $|k| \neq |m|$
$$ \int_{0}^{1/2} \widehat \chi_{\left[-s, s\right]}(k)  \widehat \chi_{\left[-s, s\right]}(m) ds = \int_{0}^{1/2}{ \frac{\sin{(2 k \pi s)}}{k \pi} \frac{\sin{(2 m \pi s)}}{m \pi}  ds} = 0.$$
Finally, we remark that
\begin{align*}
 \int_{0}^{1/2}  \widehat \chi_{\left[-s, s\right]}(k)^2 ds &=   \int_{0}^{1/2}  \widehat \chi_{\left[-s, s\right]}(k)^2 ds \\
&=\int_{0}^{1/2}{ \left( \frac{\sin{(2 k \pi s)}}{k \pi} \right)^2 ds} = \frac{1}{k^2 \pi^2} \int_{0}^{1/2}{\sin{(2 k \pi s)^2}ds} = \frac{1}{4 k^2 \pi^2}.
\end{align*}
Therefore, the expression simplifies to
$$  \frac{1}{N^2} \sum_{k \in \mathbb{Z} \atop k \neq 0}  \frac{1}{2k^2 \pi^2} \left| \sum_{n=1}^{N}{ e^{2 \pi i  k x_n }}\right|^4 =   \frac{1}{N^2} \sum_{k =1}^{\infty}  \frac{1}{k^2 \pi^2} \left| \sum_{n=1}^{N}{ e^{2 \pi i  k x_n }}\right|^4.$$
The second term simplifies to
\begin{align*}
-\frac{2}{N}  \sum_{k \in \mathbb{Z} \atop k \neq 0} \left| \sum_{n=1}^{N}{ e^{2 \pi i  k x_n }}\right|^2  \int_{0}^{1/2} \widehat \chi_{\left[-s, s\right]}(k)ds &= - \frac{2}{N}    \sum_{k \in \mathbb{Z} \setminus (2\mathbb{Z}) \atop k \neq 0}   \frac{1}{k^2 \pi^2}\left| \sum_{n=1}^{N}{ e^{2 \pi i  k x_n }}\right|^2 \\
&= - \frac{4}{N}    \sum_{k \in \mathbb{N} \setminus (2\mathbb{N})}   \frac{1}{k^2 \pi^2}\left| \sum_{n=1}^{N}{ e^{2 \pi i  k x_n }}\right|^2  
\end{align*}
and the third term is trivially $1/2$. Altogether,
\begin{align*}
A &= \frac{1}{N^2}\sum_{k =1}^{\infty}  \frac{2}{\pi^2 k^2} \left| \sum_{n=1}^{N}{ e^{2 \pi i  k x_n }}\right|^4 -    \frac{1}{N}\sum_{k \in \mathbb{N} \setminus (2\mathbb{N})}   \frac{8}{k^2 \pi^2}\left| \sum_{n=1}^{N}{ e^{2 \pi i  k x_n }}\right|^2 + 1.
\end{align*}
The argument shows that we have an essentially explicit expression for the squared deviation; the remaining difficulty is to estimate the two exponential sums.
The inequality of LeVeque \cite{leveque} bounds the second term in size from above by
$$ \frac{1}{N}\sum_{k \in \mathbb{Z} \setminus (2\mathbb{Z}) \atop k \neq 0}   \frac{2}{k^2 \pi^2}  \left| \sum_{n=1}^{N}{ e^{2 \pi i  k x_n }}\right|^2 \lesssim   
 \frac{1}{N^2} \left(N \sum_{k = 1}^{\infty}   \frac{1}{k^2 }  \left| \sum_{n=1}^{N}{ e^{2 \pi i  k x_n }}\right|^2 \right) \lesssim
N D_N^2.$$
We can use the same inequality to also conclude that
$$   \frac{1}{N^2} \sum_{k =1}^{\infty}  \frac{1}{4 \pi^2 k^2} \left| \sum_{n=1}^{N}{ e^{2 \pi i  k x_n }}\right|^4 \lesssim   \sum_{k =1}^{\infty}  \frac{1}{k^2} \left| \sum_{n=1}^{N}{ e^{2 \pi i  k x_n }}\right|^2 \lesssim N^2 D_N^2.$$
We will now compute a lower bound for this term as well: recall the Erd\H{o}s-Turan inequality
$$N D_N \leq  \inf_{K \in \mathbb{N}} \frac{N}{K+1} + 3 \sum_{k=1}^{K}\frac{1}{k}  \left| \sum_{n=1}^{N}{ e^{2 \pi i  k x_n }}\right|.$$
Denoting the right-hand side by $E_N$, we summarize the Erd\H{os}-Turan inequality as $E_N \gtrsim N D_N$ and conclude that
$$ ND_N \lesssim \sum_{k=1}^{N/E_N}\frac{1}{k}  \left| \sum_{n=1}^{N}{ e^{2 \pi i  k x_n }}\right|.$$
Using H\"older's inequality $L^{4/3} \times L^4 \rightarrow L^1$, we obtain
\begin{align*}
N D_N \lesssim \sum_{k=1}^{N/E_N}\frac{1}{k}  \left| \sum_{n=1}^{N}{ e^{2 \pi i  k x_n }}\right| &\leq \left(\sum_{k=1}^{N/E_N}\frac{1}{k^{2/3}} \right)^{3/4} \left(\sum_{k =1}^{\infty}  \frac{1}{k^2} \left| \sum_{n=1}^{N}{ e^{2 \pi i  k x_n }}\right|^4 \right)^{1/4} \\
&\lesssim  \left(\frac{N}{E_N}\right)^{\frac{1}{4}} \left(\sum_{k =1}^{\infty}  \frac{1}{k^2} \left| \sum_{n=1}^{N}{ e^{2 \pi i  k x_n }}\right|^4 \right)^{1/4}
\end{align*}
and thus
$$  \frac{1}{N^2}\sum_{k =1}^{\infty}  \frac{1}{k^2} \left| \sum_{n=1}^{N}{ e^{2 \pi i  k x_n }}\right|^4 \gtrsim N  D_N^4 E_N  \gtrsim N^2 D_N^5.$$
Finally, whenever   $ N^2 D_N^5 \gtrsim N D_N^2$ (which occurs for $D_N \gtrsim N^{-1/3}$), the positive terms dominates the negative term and
$$ A \gtrsim  \frac{1}{N^2}\sum_{k =1}^{\infty}  \frac{1}{4 \pi^2 k^2} \left| \sum_{n=1}^{N}{ e^{2 \pi i  k x_n }}\right|^4 \gtrsim N^2 D_N^5.$$
\end{proof}

\textit{Remark.} The proof has an amusing consequence if we plug in the set $\left\{0, 1/2\right\}$. We observe that
\begin{align*}
 A = 2\int_{0}^{1/2} \left( \frac{1}{N} \# \left\{ 1 \leq m \neq n \leq N: |x_m - x_n| \leq s \right\} - 2Ns \right)^2 ds =2\int_{0}^{1/2} \left( 4s\right)^2 ds = \frac{4}{3}.
\end{align*}
At the same time, the exponential sum is very easy
$$ \sum_{n=1}^{2}{e^{2\pi i k x_n}} = \begin{cases} 2 \qquad &\mbox{if}~k~\mbox{is even} \\
0 \qquad &\mbox{if}~k~\mbox{is odd} \end{cases}.$$
The formula
$$ A = \frac{1}{N^2}\sum_{k =1}^{\infty}  \frac{2}{\pi^2 k^2} \left| \sum_{n=1}^{N}{ e^{2 \pi i  k x_n }}\right|^4 -    \frac{1}{N}\sum_{k \in \mathbb{N} \setminus (2\mathbb{N})}   \frac{8}{k^2 \pi^2}\left| \sum_{n=1}^{N}{ e^{2 \pi i  k x_n }}\right|^2 + 1$$ 
obtained in the proof of Theorem 2 thus simplifies to
$$ \frac{4}{3} = \frac{1}{4} \sum_{k=1}^{\infty}{\frac{2}{\pi^2 k^2} \left(16 \cdot 1_{k~\mbox{is even}}\right)} +1 = \frac{1}{4} \sum_{k=1}^{\infty}{\frac{32}{\pi^2 (2k)^2}} + 1 = 2 \sum_{k=1}^{\infty}{\frac{1}{\pi^2 k^2}} + 1$$
and thus $\zeta(2) = \pi^2/6$.

\subsection{Proof of the Corollaries.}
The second Corollary is easy to establish. We have 
$$ 0 \leq A = \frac{1}{N^2}\sum_{k =1}^{\infty}  \frac{2}{\pi^2 k^2} \left| \sum_{n=1}^{N}{ e^{2 \pi i  k x_n }}\right|^4 -    \frac{1}{N}\sum_{k \in \mathbb{Z} \setminus (2\mathbb{Z}) \atop k \neq 0}   \frac{8}{k^2 \pi^2}\left| \sum_{n=1}^{N}{ e^{2 \pi i  k x_n }}\right|^2 + 1$$
and the desired inequality follows from a re-formulation. The first Corollary can be see as follows
\begin{align*}
A \leq \frac{1}{N^2}\sum_{k =1}^{\infty}  \frac{2}{ \pi^2 k^2} \left| \sum_{n=1}^{N}{ e^{2 \pi i  k x_n }}\right|^4 +1 =  \frac{1}{N^2}\sum_{k=1}^{\infty}  \frac{2}{ \pi^2 k^2} \left| \sum_{m, n=1}^{N}{ e^{2 \pi i  k (x_n-x_m) }}\right|^2 +1.
\end{align*}
The desired inequality then follows from the definition of diaphony $F_N$.

\section{Proof of the Proposition}
\begin{proof}[Sketch of Proof.]
We use the Erd\H{o}s-Turan inequality to conclude that 
\begin{align*}
N D_N &\lesssim  \sum_{k=1}^{N}\frac{1}{k}  \left| \sum_{n=1}^{N}{ e^{2 \pi i  k x_n }}\right| \leq  \left( \sum_{k=1}^{N}{\frac{1}{k}} \right)^{1/2} \left( \sum_{k=1}^{N} \frac{1}{k}  \left| \sum_{n=1}^{N}{ e^{2 \pi i  k x_n }}\right|^2 \right)^{1/2}\\
&\lesssim  \sqrt{\log{N}}\left(     \sum_{m,n=1}^{N}  \sum_{k=1}^{N}  \frac{\cos{(2 \pi   k (x_m-x_n) )}}{k} \right)^{1/2}
\end{align*}
The main inside is that the inner sum resembles a well-known Fourier series
$$ \sum_{k=1}^{\infty}  \frac{\cos{(2 \pi   k x)}}{k} = \log\left( \frac{1}{4 \sin^2{\pi x}}\right)$$
and there is fast convergence away from the integers. Close to the origin, we may use
$$ \left| \sum_{k=1}^{N}  \frac{\cos{(2 \pi   k x)}}{k} \right| \lesssim \log{N}$$
and the transition region can be dealt with by standard methods.
\end{proof}


\begin{thebibliography}{}

\bibitem{aist1} C. Aistleitner, G. Larcher and M. Lewko, Additive Energy and the Hausdorff dimension of the exceptional set in metric pair correlation problems. With an Appendix by Jean Bourgain. to appear in Israel J. Math

\bibitem{aist} C. Aistleitner, T. Lachmann, and F. Pausinger, Pair correlations and equidistribution,  J. Number Theory 182, 206--220 (2018).

\bibitem{cas} J. W. Cassels, A new inequality with application to the theory of diophantine approximation, Math. Ann. 126, 108--118

\bibitem{david}  H. David and H. Nagaraja, Order statistics. Third edition. Wiley Series in Probability and Statistics. Wiley-Interscience, Hoboken, NJ, 2003. 

\bibitem{grep} S. Grepstad and G. Larcher, On Pair Correlation and Discrepancy, Arch. Math. (109), pp. 143--149, 2017.


\bibitem{heath} D. R. Heath-Brown, Pair correlation for fractional parts of $\alpha n^2$, Math. Proc. Cambridge Philos. Soc. 148 (2010), 385--407.

\bibitem{kui} L. Kuipers and H. Niederreiter, 
Uniform distribution of sequences. 
Pure and Applied Mathematics. Wiley-Interscience, New York-London-Sydney, 1974. 

\bibitem{leveque} W. Leveque, 
An inequality connected with Weyl's criterion for uniform distribution. 1965 Proc. Sympos. Pure Math., Vol. VIII pp. 22--30 Amer. Math. Soc., Providence, R.I. 

\bibitem{lieb} E. Lieb, M. Loss,  Analysis. Second edition. Graduate Studies in Mathematics, 14. American Mathematical Society, Providence, RI, 2001

\bibitem{mon} H. L. Montgomery. Minimal Theta Functions. Glasgow Mathematical Journal, 30, 1988.

\bibitem{dia3} Y. Ohkubo, The diaphony of a class of infinite sequences. Probability and number theory—Kanazawa 2005, 307--322,  Adv. Stud. Pure Math., 49, Math. Soc. Japan, Tokyo, 2007. 

\bibitem{dia2} F. Pausinger and W. Schmid, 
On the asymptotics of a lower bound for the diaphony of generalized van der Corput sequences. (English summary) Monte Carlo methods and applications, 163--169, 
De Gruyter Proc. Math., De Gruyter, Berlin, 2013. 

\bibitem{dia1} F. Pillichshammer, The p-adic diaphony of the Halton sequence. Funct. Approx. Comment. Math. 49 (2013), no. 1, 91--102. 

\bibitem{pir} I. Pirsic and W. Stockinger, The Champernowne constant is not Poissonian,  arXiv:1710.09313

\bibitem{rud1} Z. Rudnick and P. Sarnak, The pair correlation function of fractional parts of
polynomials, Comm. Math. Phys. 194 (1998), 61--70.

\bibitem{rud2}   Z. Rudnick and A. Zaharescu, A metric result on the pair correlation of fractional
parts of sequences, Acta. Arith. 89 (1999), 283--293.

\bibitem{rud3} Z. Rudnick and A. Zaharescu, The distribution of spacings between fractional
parts of lacunary sequences, Forum Math. 14 (2002), 691--712.

\bibitem{stein} S. Steinerberger, Localized Quantitative Criteria for Equidistribution, Acta Arithmetica, 180 , 183--199 (2017).

\bibitem{stein1}  S. Steinerberger, Exponential Sums and Riesz energies, to appear in Journal of Number Theory

\bibitem{vdc} J. van der Corput and C. Pisot, Sur la discrepance modulo un,Indag. Math. 1, 143--153, 184--195, 260--269 (1939).

\bibitem{vin} I. Vinogradov, On the fractional parts of integral polynomials (Russian), Izu. Akad.
Nauk SSSR 20, 585--600 (1926)

\bibitem{walker} A. Walker, The Primes are not metric Poissonian, 1702.07365

\bibitem{zinterhof} P. Zinterhof, 
\"Uber einige Absch\"atzungen bei der Approximation von Funktionen mit Gleichverteilungsmethoden.  
Österreich. Akad. Wiss. Math.-Naturwiss. Kl. S.-B. II 185 (1976), no. 1-3, 121--132. 

\end{thebibliography}
\end{document}